\documentclass[14pt,a4paper]{amsart}
\usepackage[14pt]{extsizes}
\usepackage{relsize}
\usepackage{amsmath}
\usepackage{amssymb}
\usepackage{tikz-cd}
\usepackage{tikz}
\usepackage{amscd}
\usepackage[left=3cm,right=2cm,top=3cm,bottom=2cm]{geometry}
\usepackage{abstract}

\newtheorem{Theorem}{Theorem}
\newtheorem{Lemma}{Lemma}

\newtheorem{Definition}{Definition}

\newtheorem{Corollary}{Corollary}

\title{A new non-reduced moduli component of rank 2 semistable sheaves on $\mathbb{P}^{3}$}

\author[A. Lavrov]{Aleksei Lavrov}
\address{Department of Mathematics\\
National Research University
Higher School of Economics\\
6 Usacheva Street\\ 
119048 Moscow, Russia}
\email{xnikiv@gmail.com}

\begin{document}

\maketitle 

\begin{abstract}
\noindent In the present paper we describe new component of the Gieseker-Maruyama moduli space $\mathcal{M}(14)$ of coherent semistable rank-2 sheaves with Chern classes $c_1=0, \ c_2=14, \ c_3=0$ on $\mathbb{P}^{3}$ which is generically non-reduced. The construction of this component is based on the technique of elementary transformations of sheaves and the famous Mumford's example of a non-reduced component of the Hilbert scheme of smooth space curves of degree 14 and genus 24. 
\end{abstract}

\section{Introduction}

Let $\mathcal{M}(0,k,2n)$ be the Gieseker-Maruyama moduli scheme of semistable rank-2 sheaves with Chern classes $c_1=0,\ c_2=k,\ c_3=2n$ on the projective space $\mathbb{P}^3$. Denote $\mathcal{M}(k)=\mathcal{M}(0,k,0)$. By the singular locus of a given $\mathcal{O}_{\mathbb{P}^3}$-sheaf $E$ we understand the set $\mathrm{Sing}(E)=\{x\in\mathbb{P}^3\ |\ E$ is not locally free at the point $x\}$. $\mathrm{Sing}(E)$ is always a proper closed subset of $\mathbb{P}^3$ and, moreover, if $E$ is a semistable sheaf of nonzero rank, every irreducible component of $\mathrm{Sing}(E)$ has dimension at most 1. For simplicity we will not make a distinction between a stable sheaf $E$ and corresponding isomorphism class $[E]$ as a point of moduli scheme. Also by a general point of an irreducible scheme we understand a closed point belonging to some Zariski open dense subset of this scheme.

Any semistable rank-2 sheaf $[E] \in \mathcal{M}(k)$ is torsion-free, so it fits into the exact triple
\begin{equation}
0 \longrightarrow E \longrightarrow E^{\vee \vee} \longrightarrow Q \longrightarrow 0,
\end{equation}
where $E^{\vee \vee}$ is a reflexive hull of $E$ and $\text{dim}~Q \leq 1$. Conversely, take a reflexive sheaf $F$, a subscheme $X \subset \mathbb{P}^{3}$, an $\mathcal{O}_{X}$-sheaf $Q$ and a surjective morphism $\phi: F \twoheadrightarrow Q$, then one can show that the kernel sheaf $E:=\text{ker}~\phi$ is semistable when $F$ and $Q$ satisfy some mild conditions. We call the triple $(F, Q, \phi)$ an \textit{elementary trasformation} data and a sheaf $E$ an \textit{elementary transform} of $F$ along $X$. 

It is interesting that all known irreducible components of the moduli schemes $\mathcal{M}(k)$ general points of which correspond to stable sheaves with singularities were described by using elementary transformations.

More precisely, in \cite{JMT2} there were found two infinite series of irreducible components of the collection $\{ \mathcal{M}(k) \}_{k=1}^{\infty}$ which (generically) parameterize stable sheaves with singularities of dimension 0 and pure dimension 1, respectively. General points of components of the first series are elementary transforms of stable reflexive sheaves along unions of distinct points in $\mathbb{P}^{3}$, while those of the second series are elementary transforms of instanton bundles along smooth complete intersection curves.

Next, in \cite{IT} there were constructed three components of $\mathcal{M}(3)$ parameterizing sheaves with singularities of mixed dimension. General sheaves of these components are elementary transforms of stable reflexive sheaves with Chern classes $(c_2, c_3)=(2, 2), \ (2, 4)$ along a disjoint union of a projective line and a collection of points in $\mathbb{P}^{3}$. This approach was generalized in \cite{AJT} by doing elementary transformations of stable reflexive sheaves with other Chern classes along a disjoint union of a projective line and a collection of points in order to construct infinite series of components of $\mathcal{M}(-1,c_2,c_3)$.

In \cite{I} the author \footnote{the author recently changed the last name from Ivanov to Lavrov} constructed an infinite series of irreducible moduli components which includes the components parameterizing non-locally free sheaves constructed in \cite{JMT2} and \cite{IT} as special cases. General sheaves of these components are obtained by elementary transformations of stable and properly $\mu$-semistable reflexive sheaves along a sheaf $Q = L \oplus \mathcal{O}_{W}$ where $W$ is a collection of points in $\mathbb{P}^{3}$ and $L$ is a line bundle over a smooth connected curve $C$ which is either rational or a complete intersection curve.

The present paper is devoted to the description of a new component of the moduli scheme $\mathcal{M}(14)$. We prove that this component is generically non-reduced. This yields the first example of a generically non-reduced component of the Gieseker-Maruyama moduli scheme of rank 2 semistable sheaves on $\mathbb{P}^{3}$. The construction of this component essentially follows the method described in the paper \cite{I}, so ideologically the new component can be considered as a particular member of the series of components constructed in \cite{I}. More precisely, a general point of the new component is an elementary transform of the trivial sheaf along a line bundle of degree 51 over a smooth space curve of degree 14 and genus 24. These curves run through a non-reduced component of the corresponding Hilbert scheme which was for the first time described by Mumford in his classical paper \cite{M}. Roughly speaking, new component of $\mathcal{M}(14)$ is fibered over the Mumford component which leads to the non-reducedness of the new component.

The paper is organized as follows. Section 2 is devoted to the construction of a new component of $\mathcal{M}(14)$. This section mainly repeats the basic computations presented in \cite{I}. In Section 3 we remind necessary facts from the deformation theory. In Section 4 we prove that this component is generically non-reduced. More precisely, we consider the first obstruction map for the scheme $\mathcal{M}(14)$ and show that it is non-trivial (see Theorem 2 and its proof). At the end of the paper we conjecture how this result can be generalized.

\vspace{3mm}\noindent{\bf Acknowledgements.}
The work was supported in part by Young Russian Mathematics award and by the Simons Foundation. I would like to thank my supervisor A. S. Tikhomirov for usefull discussions. 

\section{Construction of a new component of $\mathcal{M}(14)$}

Let $\text{Hilb}_{14, 24}$ denote the open subset of the Hilbert scheme of $\mathbb{P}^{3}$ parameterising smooth connected curves of degree 14 and genus 24 in $\mathbb{P}^{3}$. Consider a $56$-dimensional irreducible open subscheme $\mathcal{H} \subset (\text{Hilb}_{14, 24})_{red}$ whose general member $C$ is contained in a smooth cubic surface. Mumford \cite{M} showed that the dimension of the tangent space of $\text{Hilb}_{14, 24}$ at $[C]$ is equal to 57. Moreover, he proved that $\mathcal{H}$ is maximal as a subvariety of $(\text{Hilb}_{14, 24})_{\text{red}}$, and hence $\text{Hilb}_{14, 24}$ is non-reduced at all points $\mathcal{H}$. Consider the universal curve $\mathcal{Z} \subset \mathbb{P}^{3} \times \mathcal{H}$ of the Hilbert scheme $\text{Hilb}_{14, 24}$ restricted onto the subset $\mathcal{H}$.

Note that the morphism $\mathcal{Z} \rightarrow \mathcal{H}$ is projective and flat, and its fibers are smooth and connected. So there exists a relative Picard scheme $\text{Pic}_{\mathcal{Z}/\mathcal{H}}$ with a Poincar\'e sheaf $\textbf{L}$ (see \cite[Thm. 4. 18. 1 and Ex. 4. 3]{K}). Consider a component $\widetilde{\mathcal{P}} \subset \text{Pic}_{\mathcal{Z}/\mathcal{H}}$ parameterizing invertible sheaves with the Hilber polynomial $P(k)=14(k+2)$. Finally, denote by $\mathcal{P} \subset \widetilde{\mathcal{P}}$ an open subset of $\widetilde{\mathcal{P}}$ consisting of line bundles $L$ satisfying the following equalities 
\begin{equation}\label{main_curve_eqs}
h^{1}(L) = 0, \ \ \ \ \ h^{0}(\omega_{C}(4) \otimes L^{-2}) = 0.
\end{equation}
Taking if necessary an open subset of $\mathcal{H}$ we can assume that $\mathcal{P}$ is smooth (see \cite[Prop. 5.19]{K}). We will denote also by $\textbf{L}$ the restriction of the Poincar\'e sheaf onto $\mathcal{P}$. Note that for any $(C, L) \in \mathcal{P}$ we have the equalities
\begin{equation}\label{h_L}
h^{0}(L) = \chi(L) = P(0) = 28, \ \ \ \text{deg}(L) = \chi(L) + g - 1 = 51.
\end{equation} 
Also the dimension of the scheme $\mathcal{P}$ is equal to
\begin{equation}\label{dim_P}
\text{dim}~\mathcal{P} = \text{dim}~\mathcal{H} + \text{dim Jac}(C) = 56 + 24 = 80.
\end{equation} 

Suppose that $L$ is a line bundle over a smooth curve $C$ representing some point of $\mathcal{P}$. Let $\text{Hom}_{e}\big( 2\mathcal{O}_{\mathbb{P}^{3}}, L \big) \subset \text{Hom}\big( 2\mathcal{O}_{\mathbb{P}^{3}}, L \big)$ be the open subset of surjective morphisms $2\mathcal{O}_{\mathbb{P}^{3}} \twoheadrightarrow L$. Note that the group $\text{Aut}(2\mathcal{O}_{\mathbb{P}^{3}}) \times \text{Aut}(L) \simeq GL(2) \times k^{*}$ acts on $\text{Hom}_{e}\big( 2\mathcal{O}_{\mathbb{P}^{3}}, L \big)$ by $(\psi, \zeta) \phi = \zeta \circ \phi \circ \psi^{-1}$ for $(\psi, \zeta) \in \text{Aut}(2\mathcal{O}_{\mathbb{P}^{3}}) \times \text{Aut}(L), \ \phi \in \text{Hom}_{e}\big( 2\mathcal{O}_{\mathbb{P}^{3}}, L \big)$. 

An epimorphism $\phi \in \text{Hom}_{e}\big( 2\mathcal{O}_{\mathbb{P}^{3}}, L \big)$ defines the kernel sheaf $E:=\text{ker}~\phi$. Moreover, two epimorphisms $\phi_1, \phi_2 \in \text{Hom}_{e}\big( 2\mathcal{O}_{\mathbb{P}^{3}}, L \big)$ define isomorphic kernel sheaves $\text{ker}~\phi_1 \simeq \text{ker}~\phi_2$ if and only if there exists some $(\psi, \zeta) \in \text{Aut}(2\mathcal{O}_{\mathbb{P}^{3}}) \times \text{Aut}(L)$ such that $(\psi, \zeta) \phi_1 = \phi_2$ (see \cite[Corr. 1.5]{SRS}). Therefore, an element $[\phi]$ of the orbit space 
$$
\text{Hom}_{e}\big( 2\mathcal{O}_{\mathbb{P}^{3}}, L \big) / \Big( \text{Aut}(2\mathcal{O}_{\mathbb{P}^{3}}) \times \text{Aut}(L) \Big) \simeq \mathbb{P}\text{Hom}_{e}\big( 2\mathcal{O}_{\mathbb{P}^{3}}, L \big) / PGL(2),
$$
which we will consider as a set, uniquely defines an isomorphism class $[E:=\text{ker}~\phi]$. Next, since $c_1(L)=0, \ c_2(L) = -14$ and $c_3(L) = 0$ we have $c_1(E) = 0, \ c_2(E) = 14$ and $c_3(E) = 0$. Moreover, from Lemma 4.3 in \cite{JMT1} follows that the sheaf $E$ is stable. Overall, these facts imply that the elements of the following set of data of elementary transformations
\begin{equation}
\mathcal{Q}:=\Big\{(C, L, [\phi]) \ | \ (C, L) \in \mathcal{P}, \ [\phi] \in \mathbb{P}\text{Hom}_{e}\big( 2\mathcal{O}_{\mathbb{P}^{3}}, L \big) / PGL(2) \Big\}
\end{equation}
\noindent are in one-to-one correspondence with some subset of closed points of the moduli scheme $\mathcal{M}(14)$.

\begin{Lemma}
There exists an irreducible closed subset $\overline{\mathcal{C}}$ of $\mathcal{M}(14)$ and a dense subset $\mathcal{C} \subset \overline{\mathcal{C}}$ whose closed points are in one-to-one correspondence with the elements of the set $\mathcal{Q}$. The dimension of $\overline{\mathcal{C}}$ is equal to 132.
\end{Lemma}

\noindent\textit{Proof:} Let $i: \mathcal{P} \times_{\mathcal{H}} \mathcal{Z} \hookrightarrow \mathcal{P} \times \mathbb{P}^{3}$ be the natural inclusion and $p: \mathcal{P} \times \mathbb{P}^{3} \twoheadrightarrow \mathcal{P}$ be the projection onto the first term. Consider the $\mathcal{O}_{\mathcal{P}}$-sheaf $\tau:=p_{*} \mathcal{H}om(2 \mathcal{O}_{\mathcal{P} \times \mathbb{P}^{3}}, \ i_{*} \textbf{L})$. Since the sheaf $\textbf{L}$ is flat over $\mathcal{P}$ the Euler characteristics $\chi(\textbf{L}|_{\{(C, L)\} \times C})$ is constant for all $(C, L) \in \mathcal{P}$. Moreover, by construction of $\mathcal{P}$, we have $h^{1}(L) = 0$ for any pair $(C, L) \in \mathcal{P}$, so $h^{0}(\textbf{L}|_{\{(C, L)\} \times C}) = \chi(\textbf{L}|_{\{(C, L)\} \times C})$ is also constant. This means that the fibers $\tau_{y} \otimes k(y), \ y \in \mathcal{P}$ are of the same dimension. Now taking into account that the scheme $\mathcal{P}$ is smooth we obtain that the sheaf $\tau$ is actually locally-free, so it can be considered as a vector bundle.

Consider the projective bundle $\textbf{P}(\tau^{\vee}):=\text{Proj}(\text{Sym}_{\mathcal{O}_{\mathcal{P}}}(\tau^{\vee}))$ associated to the vector bundle $\tau$. If the point $u \in \mathcal{P}$ corresponds to the pair $(C, L)$, then the fiber of the projection $\textbf{P}(\tau^{\vee}) \longrightarrow \mathcal{P}$ over the point $u$ is the projective space $\mathbb{P}\text{Hom}(2 \mathcal{O}_{\mathbb{P}^{3}},L)$. Using this observation and the formulas (\ref{dim_P}), (\ref{h_L}) we can compute dimension of $\textbf{P}(\tau^{\vee})$ as follows  
\begin{equation}\label{dim of proj}
\text{dim}~\textbf{P}(\tau^{\vee}) = \text{dim}~ \mathcal{P} + \text{dim}~ \mathbb{P}\text{Hom}(2 \mathcal{O}_{\mathbb{P}^{3}},L) = 80 + 2 \cdot 28 - 1 = 135.
\end{equation}

Let $\frak{E} \subset \textbf{P}(\tau^{\vee})$ be the open dense subset of $\textbf{P}(\tau^{\vee})$ consisting of classes of surjective morphisms $[2 \mathcal{O}_{\mathbb{P}^{3}} \twoheadrightarrow L]$. Any point $q \in \frak{E}$ determines the isomorphism class of the sheaf $[E_{q}:=\text{ker} \ \psi_{q}]$, where $[\psi_{q}] \in \mathbb{P}\text{Hom}_{e}(2 \mathcal{O}_{\mathbb{P}}, L)$. The family $\{ E_{q}, \ q \in \frak{E} \}$ globalizes to the universal sheaf $\textbf{E}$ over $\frak{E} \times \mathbb{P}^{3}$. In order to show this note that we have the anti-tautological line bundle $\mathcal{O}_{\textbf{P}(\tau^{\vee})}(1)$ over $\textbf{P}(\tau^{\vee})$ such that $\text{pr}_{*}(\mathcal{O}_{\textbf{P}(\tau^{\vee})}(1)) = \tau^{\vee}$ where $\text{pr}:\textbf{P}(\tau^{\vee}) \longrightarrow \mathcal{P}$ is the natural projection. The canonical section $\sigma \in \Gamma(\mathcal{P}, \tau \otimes \tau^{\vee})$ comes from a unique section $\widetilde{\sigma} \in \Gamma(\textbf{P}(\tau^{\vee}), \text{pr}^{*}\tau \otimes \mathcal{O}_{\textbf{P}(\tau^{\vee})}(1))$. The latter defines a morphism of sheaves $\phi: \text{pr}^{*}(2\mathcal{O}_{\mathcal{P} \times \mathbb{P}^{3}}) \longrightarrow \text{pr}^{*}i_{*}\textbf{L} \otimes \mathcal{O}_{\textbf{P}(\tau^{\vee})}(1)$ and we set $\textbf{E} := \text{ker}~\phi$. So we have the following exact triple:
$$
0 \longrightarrow \textbf{E} \longrightarrow \text{pr}^{*}(2\mathcal{O}_{\mathcal{P} \times \mathbb{P}^{3}}) \longrightarrow \text{pr}^{*}i_{*}\textbf{L} \otimes \mathcal{O}_{\textbf{P}(\tau^{\vee})}(1) \longrightarrow 0.
$$
One can see that the sheaf $\textbf{E}$ being restricted onto $\frak{E}$ is the desired globalization of the family of sheaves $\{ E_{q}, \ q \in \frak{E} \}$.

Next, by construction and by definition of moduli scheme, the sheaf $\textbf{E}$ defines the modular morphism $\Phi: \frak{E} \longrightarrow \mathcal{M}(14), \ q \mapsto [E_{q}=\text{ker} \ \psi_{q}]$. Now consider the image $\mathcal{C}:=\text{im}(\Phi)$ of the morphism $\Phi$ as a subset of $\mathcal{M}(14)$ and take its closure $\overline{\mathcal{C}} \subset \mathcal{M}(m+d)$. Note that the scheme $\frak{E}$ is irreducible, so the closed subset $\overline{\mathcal{C}}$ is also irreducible. Moreover, taking if necessary an open dense subset of $\mathcal{C}$ we can assume that the morphism $\Phi$ is flat over $\mathcal{C}$. In particular, this means that for a general point $[E] = \Phi(x), \ x \in \frak{E}$, we have the following formula for dimension
\begin{equation}\label{dim of image}
\text{dim}_{[E]}~ \overline{\mathcal{C}} = \text{dim}_{x}~\frak{E} - \text{dim}_{x}~\Phi^{-1}([E]).
\end{equation}
Moreover, we have $\Phi(x)=\Phi(y)$ if and only if the corresponding equivalence classes of morphisms $[\phi_{x}], \ [\phi_{y}] \in \mathbb{P}\text{Hom}_{e}(2 \mathcal{O}_{\mathbb{P}^{3}}, L)$ differ by the action of the group $PGL(2)$. Since this action is free, the fiber $\Phi^{-1}([E])$ is isomorphic to the group $PGL(2)$. This implies that the set of closed points of $\mathcal{C}$ is isomorphic to $\mathcal{Q}$. Moreover, from formulas (\ref{dim of proj}) and (\ref{dim of image}) it follows that the dimension of the scheme $\overline{\mathcal{C}} \subset \mathcal{M}(14)$ is equal to 132.
\hfill$\Box$

\vspace{3mm}

\begin{Theorem}\label{Main result}
For any sheaf $[E]$ from the subset $\mathcal{C} \subset \mathcal{M}(14)$ we have the equality
$$
\emph{dim}~T_{[E]}\mathcal{M}(14)=\emph{dim}~\overline{\mathcal{C}} + 1 = 133.
$$
\end{Theorem}

\noindent\textit{Proof:} For the computation of the dimension of the tangent space of the moduli scheme $\mathcal{M}(14)$ at the point $[E]$ defined above, we use the standard fact of deformation theory, $T_{[E]}\mathcal{M}(m+d) \simeq \text{Ext}^{1}(E,E)$ for a stable sheaf $E$, and the local-to-global spectral sequence $\text{H}^p(\mathcal{E}xt^{q}(E,E)) \Rightarrow \text{Ext}^{p+q}(E,E)$, which yields the following exact sequence
\begin{equation}
0 \longrightarrow \text{H}^{1}(\mathcal{H}om(E,E)) \longrightarrow \text{Ext}^{1}(E,E) \longrightarrow \text{H}^{0}(\mathcal{E}xt^{1}(E,E)) \longrightarrow
\end{equation}
$$
\overset{\phi}{\longrightarrow} \text{H}^{2}(\mathcal{H}om(E,E)) \longrightarrow \text{Ext}^{2}(E,E).
$$

According to our construction, any sheaf $[E] \in \mathcal{C}$ fits into the exact triple of the following form
\begin{equation}\label{main}
0 \longrightarrow E \longrightarrow 2 \mathcal{O}_{\mathbb{P}^{3}} \longrightarrow L \longrightarrow 0.
\end{equation}
Moreover, the first equality of (\ref{main_curve_eqs}) and the triple
\begin{equation}
\text{H}^{1}(L) \longrightarrow \text{H}^{2}(E) \longrightarrow \text{H}^{2}(2 \mathcal{O}_{\mathbb{P}^{3}})
\end{equation}
\noindent yields $\text{H}^{2}(E) = 0$. Applying the functor $\mathcal{H}om(-,E)$ to the triple (\ref{main}) we obtain the following exact sequence
\begin{equation}
0 \longrightarrow \mathcal{H}om(L, E) \longrightarrow \mathcal{H}om(2 \mathcal{O}_{\mathbb{P}^{3}}, E)  \longrightarrow \mathcal{H}om(E,E) \longrightarrow
\end{equation}
$$
\longrightarrow \mathcal{E}xt^{1}(L, E) \longrightarrow \mathcal{E}xt^{1}(2 \mathcal{O}_{\mathbb{P}^{3}}, E) \longrightarrow \mathcal{E}xt^{1}(E,E) \longrightarrow
$$
$$
\longrightarrow \mathcal{E}xt^{2}(L,E) \longrightarrow \mathcal{E}xt^{1}(2 \mathcal{O}_{\mathbb{P}^{3}}, E).
$$
Since $E$ is torsion-free sheaf we have $ \mathcal{H}om(L, E)=0$. Also we have that the sheaf $\mathcal{E}xt^{i \geq 1}(2 \mathcal{O}_{\mathbb{P}^{3}}, E)$ is zero, so we obtain the following triple and isomorphism
\begin{equation}\label{triple2}
0 \longrightarrow \mathcal{H}om(2 \mathcal{O}_{\mathbb{P}^{3}}, E)  \longrightarrow \mathcal{H}om(E, E) \longrightarrow \mathcal{E}xt^{1}(L, E) \longrightarrow 0,
\end{equation}
\begin{equation}\label{isom1 for ext1}
\mathcal{E}xt^{1}(E,E) \simeq \mathcal{E}xt^{2}(L,E).
\end{equation}
Taking into account that the sheaf $\mathcal{E}xt^{1}(L, E)$ has dimension at most $1$ and $h^{2}(E) = 0$, we have that $h^{2}(\mathcal{H}om(E, E) = 0$. Therefore, we obtain the following exact triple
\begin{equation}\label{group_ext_1st}
0 \longrightarrow \text{H}^{1}(\mathcal{H}om(E,E)) \longrightarrow \text{Ext}^{1}(E,E) \longrightarrow \text{H}^{0}(\mathcal{E}xt^{1}(E,E)) \longrightarrow 0.
\end{equation}

Next, apply the functor $\mathcal{H}om(L,-)$ to the triple (\ref{main})
\begin{equation}
0 \longrightarrow \mathcal{H}om(L, E) \longrightarrow \mathcal{H}om(L, 2 \mathcal{O}_{\mathbb{P}^{3}}) \longrightarrow \mathcal{H}om(L, L) \longrightarrow
\end{equation}
$$
\longrightarrow \mathcal{E}xt^{1}(L, E) \longrightarrow \mathcal{E}xt^{1}(L, 2 \mathcal{O}_{\mathbb{P}^{3}}) \longrightarrow \mathcal{E}xt^{1}(L, L) \longrightarrow
$$
$$
\longrightarrow \mathcal{E}xt^{2}(L, E) \longrightarrow \mathcal{E}xt^{2}(L, 2 \mathcal{O}_{\mathbb{P}^{3}}) \longrightarrow \mathcal{E}xt^{2}(L ,L) \longrightarrow 0.
$$
Since the sheaves $E$ and $2 \mathcal{O}_{\mathbb{P}^{3}}$ are torsion-free we have that $\mathcal{H}om(L, E) = \mathcal{H}om(L, 2 \mathcal{O}_{\mathbb{P}^{3}}) = 0$. Note that the smooth curve $C$ is locally complete intersection, so for any point $x \in \mathbb{P}^{3}$ we have the following
\begin{equation}\label{local_isom}
\text{Ext}^1_{\mathcal{O}_{\mathbb{P}^{3},x}}(\mathcal{O}_{C,x},\mathcal{O}_{\mathbb{P}^{3},x}) = 0.
\end{equation}
This equality implies that the sheaf $\mathcal{E}xt^{1}(L, 2 \mathcal{O}_{\mathbb{P}^{3}})$ is equal to zero, so we have the isomorphism
\begin{equation}\label{isom1}
\mathcal{E}xt^{1}(L, E) \simeq \mathcal{H}om(L, L).
\end{equation}
and the exact sequence
\begin{equation}\label{exact seq for ext2_0}
0 \longrightarrow \mathcal{E}xt^{1}(L, L) \longrightarrow \mathcal{E}xt^{2}(L, E) \longrightarrow \mathcal{E}xt^{2}(L, 2 \mathcal{O}_{\mathbb{P}^{3}}) \longrightarrow \mathcal{E}xt^{2}(L,L) \longrightarrow 0.
\end{equation}
The following exact triple holds
\begin{equation}\label{restriction of refl}
0 \longrightarrow L^{-1} \longrightarrow 2 \mathcal{O}_{\mathbb{P}^{3}} \otimes \mathcal{O}_{C} \longrightarrow L \longrightarrow 0.
\end{equation}
Note that $\mathcal{E}xt^{2}(L, 2 \mathcal{O}_{\mathbb{P}^{3}}) \simeq \mathcal{E}xt^{2}(L, 2 \mathcal{O}_{C})$. In particular, it means that $\mathcal{E}xt^{2}(L, 2 \mathcal{O}_{\mathbb{P}^{3}})$ is locally-free $\mathcal{O}_{C}$-sheaf. Since $\mathcal{E}xt^{2}(L,L^{-1})$ is also locally-free $\mathcal{O}_{C}$-sheaf of rank 1, then applying the functor $\mathcal{H}om(L,-)$ to the triple (\ref{restriction of refl}) we obtain the following exact triple
\begin{equation}\label{ex_tr_1}
0 \longrightarrow \mathcal{E}xt^{2}(L, L^{-1}) \longrightarrow \mathcal{E}xt^{2}(L, 2 \mathcal{O}_{\mathbb{P}^{3}}) \longrightarrow \mathcal{E}xt^{2}(L, L) \longrightarrow 0.
\end{equation}
Moreover, we have the following commutative diagram
\begin{equation*}
\begin{tikzcd}[column sep=small]
\mathcal{E}xt^{2}(L, 2 \mathcal{O}_{\mathbb{P}^{3}})  \rar \dar{\simeq} & \mathcal{E}xt^{2}(L, L) \dar{=} \\
\mathcal{E}xt^{2}(L, 2 \mathcal{O}_{\mathbb{P}^{3}} \otimes \mathcal{O}_{C})  \rar & \mathcal{E}xt^{2}(L, L)
\end{tikzcd}
\end{equation*}
So the morphism $\mathcal{E}xt^{2}(L, 2 \mathcal{O}_{\mathbb{P}^{3}}) \longrightarrow \mathcal{E}xt^{2}(L,L)$ in the triple (\ref{ex_tr_1}) coincides with the last morphism in the exact sequence (\ref{exact seq for ext2_0}). Therefore, we can simplify (\ref{exact seq for ext2_0}) as 
\begin{equation}\label{exact seq for ext2}
0 \longrightarrow \mathcal{E}xt^{1}(L, L) \longrightarrow \mathcal{E}xt^{2}(L, E) \longrightarrow \mathcal{E}xt^{2}(L, L^{-1}) \longrightarrow 0.
\end{equation}
Note that for any subscheme $C \subset \mathbb{P}^{3}$ we have that 
$$
\mathcal{E}xt^{1}(\mathcal{O}_{C}, \mathcal{O}_{C}) \simeq \mathcal{H}om(I_{C}, \mathcal{O}_{C}) \simeq \mathcal{H}om(I_{C}/I_{C}^{2}, \mathcal{O}_{C}) = N_{C/\mathbb{P}^{3}}.
$$ 
(the last equality is the definition of the normal sheaf). Besides, if $C$ is a locally complete intersection of the pure dimension $1$ then (see \cite[Prop. 7.5]{AG})
$$
\mathcal{E}xt^{2}(\mathcal{O}_{C}, \mathcal{O}_{C}) \simeq \mathcal{E}xt^{2}(\mathcal{O}_{C}, \mathcal{O}_{\mathbb{P}^{3}}) \simeq \mathcal{E}xt^{2}(\mathcal{O}_{C}, \omega_{\mathbb{P}^{3}})(4) \simeq \omega_{C}(4).
$$
Since $L$ is an invertible $\mathcal{O}_{C}$-sheaf it follows that 
$$
\mathcal{E}xt^{1}(L,L) \simeq \mathcal{E}xt^{1}(\mathcal{O}_{C},\mathcal{O}_{C}), \ \ \ \mathcal{E}xt^{2}(L,L^{-1}) \simeq \mathcal{E}xt^{2}(\mathcal{O}_{C}, \mathcal{O}_{C}) \otimes L^{-2}.
$$
From these formulas one can deduce the isomorphisms
$$
\mathcal{E}xt^{1}(L, L) \simeq N_{C / \mathbb{P}^{3}}, \ \ \ \mathcal{E}xt^{2}(L,L^{-1}) \simeq \omega_{C}(4) \otimes L^{-2}.
$$
Substituting them to (\ref{exact seq for ext2}) we obtain the following exact triple
\begin{equation}\label{triple1 for ext1}
0 \longrightarrow N_{C / \mathbb{P}^{3}} \longrightarrow \mathcal{E}xt^{2}(L,E) \longrightarrow \omega_{C}(4) \otimes L^{-2} \longrightarrow 0.
\end{equation}
Therefore, using the isomorphism (\ref{isom1 for ext1}), the triple (\ref{triple1 for ext1}) and the condition $h^0(\omega_{C}(4) \otimes L^{-2}) = 0$ from (\ref{main_curve_eqs}) we obtain the isomorphism
\begin{equation}\label{h0}
\text{H}^0(\mathcal{E}xt^{1}(E,E)) \simeq \text{H}^{0}(N_{C/\mathbb{P}^{3}}) \simeq T_{C}\text{Hilb}_{14,24}.
\end{equation}

After applying the functor $\mathcal{H}om(-, 2 \mathcal{O}_{\mathbb{P}^{3}})$ to the exact triple (\ref{main}) we obtain the long exact sequence of sheaves
\begin{equation}\label{v2}
0 \longrightarrow \mathcal{H}om(L, 2 \mathcal{O}_{\mathbb{P}^{3}}) \longrightarrow \mathcal{H}om(E, 2 \mathcal{O}_{\mathbb{P}^{3}}) \longrightarrow \mathcal{H}om(2 \mathcal{O}_{\mathbb{P}^{3}}, 2 \mathcal{O}_{\mathbb{P}^{3}}) \longrightarrow
\end{equation}
$$
\longrightarrow \mathcal{E}xt^{1}(L, 2 \mathcal{O}_{\mathbb{P}^{3}}) \longrightarrow \mathcal{E}xt^{1}(E, 2 \mathcal{O}_{\mathbb{P}^{3}}) \longrightarrow 0 \longrightarrow
$$
$$
\longrightarrow \mathcal{E}xt^{2}(L, 2 \mathcal{O}_{\mathbb{P}^{3}}) \longrightarrow \mathcal{E}xt^{2}(E, 2 \mathcal{O}_{\mathbb{P}^{3}}) \longrightarrow 0.
$$
As it was already explained $\mathcal{H}om(L, 2 \mathcal{O}_{\mathbb{P}^{3}})=\mathcal{E}xt^{1}(L, 2 \mathcal{O}_{\mathbb{P}^{3}})=0$ and the sheaf $\mathcal{E}xt^{1}(2 \mathcal{O}_{\mathbb{P}^{3}}, 2 \mathcal{O}_{\mathbb{P}^{3}})$ is zero, so we have the following isomorphisms
\begin{equation}\label{isom2}
\mathcal{H}om(E, 2 \mathcal{O}_{\mathbb{P}^{3}}) \simeq \mathcal{H}om(2 \mathcal{O}_{\mathbb{P}^{3}}, 2 \mathcal{O}_{\mathbb{P}^{3}}), \ \ \ \mathcal{E}xt^{1}(E, 2 \mathcal{O}_{\mathbb{P}^{3}}) \simeq 0.
\end{equation}

Consider the part of the commutative diagram with exact rows and columns obtained by applying the bifunctor $\mathcal{H}om(-,-)$ and its derivative $\mathcal{E}xt(-,-)$ to the exact triple (\ref{main}) which looks as follows
\begin{equation*}
\begin{tikzcd}[column sep=small]
  & 0  \dar & 0 \dar &  & & \\
0 \rar & \mathcal{H}om(2 \mathcal{O}_{\mathbb{P}^{3}}, E) \rar \dar & \mathcal{H}om(2 \mathcal{O}_{\mathbb{P}^{3}}, 2 \mathcal{O}_{\mathbb{P}^{3}}) \rar \dar & \mathcal{H}om(2 \mathcal{O}_{\mathbb{P}^{3}}, L) \rar & 0 & \\
0 \rar & \mathcal{H}om(E, E) \rar{\tau} \dar & \mathcal{H}om(E, 2 \mathcal{O}_{\mathbb{P}^{3}}) \dar & & & \\
 & \mathcal{E}xt^{1}(L, E) \dar  & 0  & & & \\
  & 0  &  &
\end{tikzcd}
\end{equation*}
Due to the isomorphism (\ref{isom2}) the sheaf $\text{coker}~\tau$ fits into the exact triple
\begin{equation}\label{coker1}
0 \longrightarrow \mathcal{H}om(E, E) \longrightarrow \mathcal{H}om(2 \mathcal{O}_{\mathbb{P}^{3}}, 2 \mathcal{O}_{\mathbb{P}^{3}}) \longrightarrow \text{coker}~\tau \longrightarrow 0.
\end{equation}
So we have the following exact sequence
\begin{equation}
0 \longrightarrow \text{End}(E) \longrightarrow \text{End}(2 \mathcal{O}_{\mathbb{P}^{3}}) \longrightarrow \text{H}^{0}(\text{coker}~\tau) \longrightarrow \text{H}^{1}(\mathcal{H}om(E, E)) \longrightarrow 0.
\end{equation}
On the other hand, applying the Snake Lemma to the commutative diagram above and using the isomorphism (\ref{isom1}) we have the exact triple
\begin{equation}\label{coker2}
0 \longrightarrow \mathcal{H}om(L, L) \longrightarrow \mathcal{H}om(2 \mathcal{O}_{\mathbb{P}^{3}}, L) \longrightarrow \text{coker}~\tau \longrightarrow 0.
\end{equation}
By the definition of the subscheme $\mathcal{P} \subset \text{Pic}^{51}_{\mathcal{Z}/\mathcal{H}}$ we have $h^{1}(L)=0$, so $h^1(\mathcal{H}om(2 \mathcal{O}_{\mathbb{P}^{3}}, L)) = 0$. Therefore, from the triple (\ref{coker2}) we obtain the following exact sequence
\begin{equation}\label{coker2}
0 \longrightarrow \text{End}(L) \longrightarrow \text{Hom}(2 \mathcal{O}_{\mathbb{P}^{3}}, L) \longrightarrow \text{H}^{0}(\text{coker}~\tau) \longrightarrow \text{H}^{1}(\mathcal{O}_{C}) \longrightarrow 0.
\end{equation}
Using these equalities and the fact that the sheaf $E$ is simple due to its stability, the triple (\ref{coker1}) implies the following isomorphism
\begin{equation}\label{h1}
\text{H}^{1}(\mathcal{H}om(E, E)) \simeq \text{H}^{1}(\mathcal{O}_{C}) \oplus
\end{equation}
$$
\oplus \Bigg( \bigg(\text{Hom}(2 \mathcal{O}_{\mathbb{P}^{3}}, L) / \text{End}(L) \bigg) / \bigg( \text{End}(2 \mathcal{O}_{\mathbb{P}^{3}}) / \text{End}(E) \bigg) \Bigg).
$$
Substituting the isomorphisms (\ref{h0}) and (\ref{h1}) to the exact triple (\ref{group_ext_1st}) we obtain a non canonical isomorphism
\begin{equation}\label{main_isom}
\text{Ext}^1(E,E) \simeq T_{C}\text{Hilb}_{14, 24} \oplus \text{H}^{1}(\mathcal{O}_{C}) \oplus
\end{equation}
$$
\oplus \Bigg( \bigg(\text{Hom}(2 \mathcal{O}_{\mathbb{P}^{3}}, L) / \text{End}(L) \bigg) / \bigg( \text{End}(2 \mathcal{O}_{\mathbb{P}^{3}}) / \text{End}(E) \bigg) \Bigg).
$$
Since the sheaf $E$ is stable it is also simple, so we have $\text{dim End}(E) = 1$. Therefore, the dimension of $\text{dim Ext}^{1}(E, E)$ is equal to 
\begin{equation}
\text{dim Ext}^{1}(E, E) = \text{dim }T_{C}\text{Hilb}_{14, 24} + \text{dim H}^{1}(\mathcal{O}_{C}) + 
\end{equation}
$$
+ \bigg( 2 h^{0}(L) - 1 \bigg) - \bigg( 4 - 1 \bigg) = 57 + 24 + 55 - 3 = 133. 
$$
So we have the statement of the Theorem
\begin{equation}
\text{dim}~T_{[E]}\mathcal{M}(14)=\text{dim}~\overline{\mathcal{C}} + 1.
\end{equation}
\hfill$\Box$

\section{Necessary facts from deformation theory}

Let us recall some technical details about Ext-groups. Suppose we are given by two objects $A, \ B$ from an abelian category $\textbf{C}$. Recall that elements of the group $\text{Ext}^{n}_{\textbf{C}}(A, B)$ can be represented by equivalence classes of $n$-extensions. More precisely, $\text{Ext}^{0}(A, B) = \text{Hom}(A, B)$. Next, $\text{Ext}^{1}(A, B)$ is the set of equivalence classes of extensions of $A$ by $B$, forming an abelian group under the Baer sum. Finally, the elements of higher Ext-groups $\text{Ext}^{n}(A, B)$ can be defined as equivalence classes of $n$-extensions, which are exact sequences
$$
0 \to B \to X_{n} \to \cdots \to X_{1} \to A \to 0,
$$
under the equivalence relation generated by the relation that identifies two extensions
$$
\xi : 0 \to B \to X_{n} \to \cdots \to X_{1} \to A \to 0,
$$
$$
\xi' : 0 \to B \to X'_{n} \to \cdots \to X'_{1} \to A \to 0,
$$
if there are maps $X_{m} \to X'_{m}$ for all $m$ in $\{ 1, 2, ..., n \}$ so that every resulting square commutes, that is, if there is a chain map $\xi \to \xi'$ which is the identity on $A$ and $B$.

Next, we have the following pairing between Ext-groups which is called \textit{Yoneda product} and denoted by $-\cup-$
$$
\text{Ext} ^{m}(L,M) \otimes \text{Ext}^{n}(M,N) \to \text{Ext}^{n+m}(L,N).
$$
This pairing is induced by the map
$$
\text{Hom}(L, M) \otimes \text{Hom}(M, N) \rightarrow \text{Hom}(L, N), \ \ f \otimes g \mapsto g \circ f.
$$
In terms of extensions, it can be described as follows. Suppose that an element $\alpha \in \text{Ext}^{n}(L, M)$ is represented by the extension
$$
\alpha :0 \rightarrow M \rightarrow E_{0} \rightarrow \cdots \rightarrow E_{n-1} \rightarrow L \rightarrow 0,
$$
and an element $\beta \in \text{Ext}^{m}(L, M)$ by the extension
$$
\beta : 0 \rightarrow N \rightarrow F_{0} \rightarrow \cdots \rightarrow F_{m-1} \rightarrow M \rightarrow 0.
$$
Then the Yoneda product $\alpha \cup \beta \in \text{Ext}^{m+n}(L,N)$ is represented by the concatenated extension
$$
\alpha \cup \beta :0\rightarrow N\rightarrow F_{0}\rightarrow \cdots \rightarrow F_{m-1} \rightarrow  E_{0}\rightarrow \cdots \rightarrow E_{n-1}\rightarrow L\rightarrow 0.
$$

Now move on to some basic facts on the infinitesimal study of the Hilbert scheme of space curves. Let $C$ be a smooth connected curve in $\mathbb{P}^{3}$ whose structure sheaf $\mathcal{O}_{C}$ is defined by the exact triple
\begin{equation}\label{structure_sheaf}
0 \longrightarrow I_{C} \longrightarrow \mathcal{O}_{\mathbb{P}^{3}} \longrightarrow \mathcal{O}_{C} \longrightarrow 0.
\end{equation}
Then an (embedded) $n$-th order (infinitesimal) deformation of $C \subset \mathbb{P}^{3}$ is a closed subscheme $\mathcal{C}_{n}$ of $\mathbb{P}^{3} \times \text{Spec}(k[t]/(t^{n+1}))$ which is flat over $k[t]/(t^{n+1})$ and $\mathcal{C}_{n} \otimes_{k[t]/(t^{n+1})} k = C$. The set of all first order deformations of $C \subset \mathbb{P}^{3}$ is the Zariski tangent space of of the Hilbert scheme at the point $[C]$. Let $\mathcal{C}_{1}$ be a first order deformation of $C \subset \mathbb{P}^{3}$. If there exists no second order deformation $\mathcal{C}_{2}$ of $C \subset \mathbb{P}^{3}$ such that $\mathcal{C}_{2} \otimes_{k[t]/(t^3)} k[t]/(t^2) = \mathcal{C}_{1}$, we say $\mathcal{C}_{1}$ is \textit{obstructed at the second order}. As we know the set of all first order deformations of $C \subset \mathbb{P}^{3}$ is parametrized by $\text{Hom}(I_{C}, \mathcal{O}_{C})$. Applying the functor $\text{Hom}(I_{C}, - )$ to the triple (\ref{structure_sheaf}) we obtain the isomorphism
$$
\delta : \text{Hom}(I_{C}, \mathcal{O}_{C}) \overset{\simeq}{\longrightarrow} \text{Ext}^{1}(I_{C}, I_{C}).
$$
Let $\phi \in \text{Hom}(I_{C}, \mathcal{O}_{C})$ be a first order deformation of $C \subset \mathbb{P}^{3}$. Then $\phi$ is obstructed at the second order if and only if the Yoneda product $o(\phi) := \delta(\phi) \cup_{1} \phi$ by
\begin{equation}\label{cup_Mumford}
\cup_{1}: \text{Ext}^{1}(I_{C}, I_{C}) \times \text{Hom}(I_{C}, \mathcal{O}_{C}) \longrightarrow \text{Ext}^{1}(I_{C}, \mathcal{O}_{C})
\end{equation}
\noindent is non-zero. An element $o(\phi) \in \text{Ext}^{1}(I_{C}, \mathcal{O}_{C})$ is called an obstruction to extend $\phi$ to second order deformations. 

Finally, consider the deformation theory of coherent sheaves. Let $X$ be a scheme over $k$, $F$ be a coherent sheaf on $X$ and $S$ be a some scheme over $k$ with a marked point $0$. Introduce the following definition

\begin{Definition}\label{Deformation}
Deformation of the sheaf $F$ over $S$ is a sheaf $\textbf{F}$ over $X \times \text{D}$ such that $\textbf{F}$ is flat over $D$, equipped with a homomorphism $p: \textbf{F} \rightarrow i_{*}F$ which induces an isomorphism $\textbf{F} \otimes i_{S*}\mathcal{O}_{X} \simeq i_{S*}F$ where $i_{S} : X \simeq X \times \{0\} \hookrightarrow X \times S$.
\end{Definition}

Denote by $D := \text{Spec}(k[t]/t^{2})$ the spectrum of the ring of dual numbers. Note that there is the one-to-one correspondence between the elements of the group $\text{Ext}^{1}(F, F)$ and the pairs $(\textbf{F}, p)$ where $\textbf{F}$ is a deformation of the sheaf $F$ over $D$ and $p: \textbf{F} \rightarrow i_{*}F$ is a morphism which induces an isomorphism $\textbf{F} \otimes i_{*}\mathcal{O}_{X} \simeq i_{*}F$. In particular, if we denote by $\text{pr}:X \times D \rightarrow X$ the projection onto the first term then we obtain the exact triple
\begin{equation}\label{extensions}
0 \longrightarrow F \longrightarrow \text{pr}_{*}\textbf{F} \longrightarrow F \longrightarrow 0,
\end{equation}
which can be considered as an element of the group $\text{Ext}^{1}(F, F)$. 

%\begin{Definition}
%Let $\textbf{F}_1$ and $\textbf{F}_2$ be some deformations over $S_1$ and $S_2$, respectively. A morphism between deformations $\textbf{F}_1$ and $\textbf{F}_2$ is a morphism $%%\phi: S_1 \rightarrow S_2$ such that we have  $\phi^{*} \textbf{F}_2 \simeq \textbf{F}_1$ 
%\end{Definition}

Next, let $(S, 0)$ be a germ of scheme $S$ at a point $0 \in S$. Deformation of a sheaf over a germ of scheme can be defined in the similar way as in Definition \ref{Deformation}. Recall that \textit{a morphism} between deformations $\textbf{F}_1$ and $\textbf{F}_2$ of a sheaf $F$ over $(S_1, 0)$ and $(S_2,0)$, respectively, can be defined as a morphism of germs $\phi : (S_1, 0) \rightarrow (S_2, 0)$ such that there exists an isomorphism  $\phi^{*}\textbf{F}_2 \simeq \textbf{F}_1$ satisfying the following commutative diagram
$$
\begin{tikzcd}
\phi^{*} \textbf{F}_2 \dar \rar{\simeq} & \textbf{F}_1 \dar \\
\phi^{*}i_{S_{2}*}F  \rar[equal] & i_{S_{1}*}F 
\end{tikzcd} \ \ \ \ \ \
$$

\begin{Definition}
Deformation $\textbf{F}$ of the sheaf $F$ over the germ of scheme $(S, 0)$ is called versal if it satisfies the following property. Let $\textbf{F}'$ and $\textbf{F}''$ be some deformations over $(S', 0)$ and $(S'',0)$, respectively. Suppose that we have a morphism $\phi:(S',0) \rightarrow (S'',0)$ from $\textbf{F}'$ to $\textbf{F}''$ and a $\psi: (S',0) \rightarrow (S,0)$ from $\textbf{F}'$ to $\textbf{F}$. Then there exists a morphism $\theta: (S'',0) \rightarrow (S,0)$ from $\textbf{F}''$ to $\textbf{F}$ such that $\psi = \theta \circ \phi$.
\end{Definition}

To any deformation $\textbf{F}$ of a sheaf $F$ over $(S, 0)$ we can associate the Kodaira-Spencer map $\tau_{\textbf{F}}: T_{0} S \rightarrow \text{Ext}^{1}(F, F)$ from the tangent space $T_{0} S$ to the set of deformations of $F$ over $D$ which can be identified with the group $\text{Ext}^{1}(F, F)$ as it was mentioned above.

\begin{Definition}
Deformation $\textbf{F}$ of the sheaf $F$ over the germ of scheme $(S, 0)$ is called semi-universal if it is versal and the corresponding Kodaira-Spencer map $\tau_{\textbf{F}}$ is an isomorphism.
\end{Definition}

\begin{Lemma}\label{obstruction_map} Let $X$ be a smooth projective variety, $F$ a coherent sheaf on $X$. Then there exists a germ of a nonsingular algebraic variety $(M, 0)$ together with a morphism $\mathcal{Y}_{F}:(M, 0) \rightarrow (\emph{Ext}^{2}(F, F), 0)$, called the obstruction map, such that the following properties are verified:

\noindent (i) $(\mathcal{Y}^{-1}_{F}(0), 0)$ is the base of a semi-universal deformation $\textbf{F}$ of the sheaf $F$. The Kodaira-Spencer map of this deformation provides a natural isomorphism $\text{KS} : T_{0}M \xrightarrow{\sim} \emph{Ext}^{1}(F, F)$.

\noindent (ii) Let 
$$
\mathcal{Y}_{F} = \sum\limits_{i=0}^{\infty}\mathcal{Y}_{F, i}, \ \ \mathcal{Y}_{F, i} \in \emph{Hom}\bigg( \text{S}^{i}(T_{0}M), \ \emph{Ext}^{2}(F, F) \bigg)
$$ 
be a Taylor expansion of $\mathcal{Y}_{F}$. Then $\mathcal{Y}_{F, 1} = 0$ and $\mathcal{Y}_{F, 2}$ is the composition
$$
T_{0}M \overset{(\text{KS}, \ \text{KS})}{\xrightarrow{\hspace*{1.6cm}}} \emph{Ext}^{1}(F, F) \times \emph{Ext}^{1}(F, F) \overset{(\xi, \eta) \mapsto \xi \cup \eta}{\xrightarrow{\hspace*{1.6cm}}} \emph{Ext}^{2}(F, F) 
$$
where $\xi \cup \eta$ denotes the Yoneda product of two elements of $\emph{Ext}^{1}(F, F)$.
\end{Lemma}
\noindent\textit{Proof:} The Appendix of Bingener to \cite{BH} provides the following scheme of the proof of this statement. The existence of a formal versal deformation was proven in \cite{Rim}. By \cite{Art}, the formal versal deformation is the formal completion of a genuine versal deformation. For the construction of $\mathcal{Y}_{F, i}$ for all $i$, see Proposition A.1 of \cite{LS}. The identification of the obstruction $\mathcal{Y}_{F, 2}$ on the formal level with the Yoneda product was done in \cite{Ar, Mu-2}.
\hfill$\Box$

\vspace{3mm}

\section{The proof of non-reducibility}

Below we will prove that the obstruction map $\mathcal{Y}_{E}$ for any sheaf $[E] \in \mathcal{C}$ is not trivial. This will immediatelly imply the main result of the paper. The idea of the proof is to connect the obstruction map $\mathcal{Y}_{E}$ for the sheaf $E$ with the obstruction map $\mathcal{Y}_{L}$ for the sheaf $L$ from the triple (\ref{main}). First of all, let us make sure that the obstruction map $\mathcal{Y}_{L}$ is not trivial.

\begin{Lemma}\label{Mumford_obstruction}
	For any sheaf $L \in \mathcal{P}$ the following Yoneda product
	\begin{equation}\label{L_cup}
	\emph{Ext}^{1}(L, L) \times \emph{Ext}^{1}(L, L) \overset{\cup_{2}}{\longrightarrow} \emph{Ext}^{2}(L, L)
	\end{equation}
	is non-trivial. Therefore, the first obstruction map $\mathcal{Y}_{L, 2}$ is non-zero.
\end{Lemma}
\noindent\textit{Proof}: According to \cite[Prop. 3. 1, remark on the page 133]{N}, any curve $C \in \mathcal{H}$ has an embedded first order infinitesimal deformation $\phi_{C} \in \text{Hom}(I_{C}, \mathcal{O}_{C})$ which is obstructed at the second order. In terms of the previous section this means that the Yoneda product (\ref{cup_Mumford}) is not trivial. 

On the other hand, we have the isomorphism $\xi_{1, 2} : \text{Ext}^{1, 2}(\mathcal{O}_{C}, \mathcal{O}_{C}) \simeq \text{Ext}^{1, 2}(L, L)$. Note that from the local-to-global spectral sequence it follows that
$$
\text{dim Ext}^{1}(\mathcal{O}_{C}, \mathcal{O}_{C}) = h^{1}(\mathcal{H}om(\mathcal{O}_{C}, \mathcal{O}_{C})) + h^{0}(\mathcal{E}xt^{1}(\mathcal{O}_{C}, \mathcal{O}_{C})).
$$
Since $\mathcal{H}om(\mathcal{O}_{C}, \mathcal{O}_{C}) \simeq \mathcal{O}_{C}$ and $\mathcal{E}xt^{1}(\mathcal{O}_{C}, \mathcal{O}_{C}) \simeq \mathcal{H}om(I_{C}, \mathcal{O}_{C})$ we have
\begin{equation}\label{lm_dim_ext}
\text{dim Ext}^{1}(\mathcal{O}_{C}, \mathcal{O}_{C}) = h^{1}(\mathcal{O}_{C}) + \text{dim Hom}(I_{C}, \mathcal{O}_{C}).
\end{equation}
Now applying the functor $\text{Hom}(-, \mathcal{O}_{C})$ to the triple (\ref{structure_sheaf}) and using the equality (\ref{lm_dim_ext}) we obtain the following exact triple
$$
0 \longrightarrow \text{Hom}(I_{C}, \mathcal{O}_{C}) \overset{\sigma}{\longrightarrow} \text{Ext}^{1}(\mathcal{O}_{C}, \mathcal{O}_{C}) \longrightarrow \text{H}^{1}(\mathcal{O}_{C}) \longrightarrow 0, 
$$
and the isomorphism $\zeta : \text{Ext}^{1}(I_{C}, \mathcal{O}_{C}) \simeq \text{Ext}^{2}(\mathcal{O}_{C}, \mathcal{O}_{C})$. Using the maps $\delta, \ \xi_{1, 2}, \ \sigma, \ \zeta$ we obtain the following commutative diagram
\begin{equation}
\begin{tikzcd}
\text{Ext}^{1}(I_{C}, I_{C}) \times \text{Hom}(I_{C}, \mathcal{O}_{C}) \rar{\cup_{1}} \dar{(\xi_{1} \circ \sigma \circ \delta^{-1}) \times (\xi_{1} \circ \sigma)} & \text{Ext}^{1}(I_{C}, \mathcal{O}_{C}) \dar{\xi_{2} \circ \zeta} \\
\text{Ext}^{1}(L, L) \times \text{Ext}^{1}(L, L) \rar{\cup_{2}} & \text{Ext}^{2}(L, L)
\end{tikzcd}
\end{equation} 
Since the left vertical map is injective, the right vertical map is an isomorphism and the upper horizontal map is non-trivial we obtain that the bottom horizontal map is also non-trivial.
\hfill$\Box$

Now in order to establish a connection between the obstruction map $\mathcal{Y}_{E}$ for the sheaf $E$ with the obstruction map $\mathcal{Y}_{L}$ for the sheaf $L$ we need the following technical lemma.

\begin{Lemma}\label{square_sequence}
Suppose that we have the following commutative diagram
\begin{equation}\label{diagram_over_dual}
\begin{tikzcd}
 & 0 \dar & 0 \dar & 0 \dar & \\
0 \rar & E \rar{\xi} \dar{\phi_1} & F \rar{\zeta} \dar{\phi_2} & L \rar \dar{\phi_3} & 0 \\
0 \rar & E_{0} \rar{\eta} \dar{{\psi_1}} & F_{0} \rar{\gamma} \dar{\psi_2} & L_{0} \rar \dar{\psi_3} & 0 \\
0 \rar & E \rar{\xi} \dar & F \rar{\zeta} \dar & L \rar \dar & 0 \\
 & 0 & 0 & 0 & 
\end{tikzcd}
\end{equation}
Each vertical extension in this diagram represents the element of the group $\emph{Ext}^{1}(E,E)$, $\emph{Ext}^{1}(F,F)$, and $\emph{Ext}^{1}(L,L)$, respectively. Consider these extensions as 3-complexes and denote them by $\mathcal{E}, \ \mathcal{F}, \ \mathcal{L}$. Then the Yoneda squares of these complexes also fits into the exact triple
\begin{equation}\label{Yoneda_square_formula}
0 \longrightarrow  \mathcal{E} \cup \mathcal{E}  \longrightarrow \mathcal{F} \cup \mathcal{F} \longrightarrow \mathcal{L} \cup \mathcal{L} \longrightarrow 0.
\end{equation}
\end{Lemma}
\noindent\textit{Proof:} The triple $(\ref{Yoneda_square_formula})$ of 4-complexes can be read as the following diagram of sheaves
\begin{equation}
\begin{tikzcd}
 & 0 \dar & 0 \dar & 0 \dar & \\
0 \rar & E \rar{\xi} \dar{\phi_1} & F \rar{\zeta} \dar{\phi_2} & L \rar \dar{\phi_3} & 0 \\
0 \rar & E_{0} \rar{\eta} \dar{\phi_1 \circ \psi_1} & F_{0} \rar{\gamma} \dar{\phi_2 \circ \psi_2} & L_{0} \rar \dar{\phi_3 \circ \psi_3} & 0 \\
0 \rar & E_{0} \rar{\eta} \dar{{\psi_1}} & F_{0} \rar{\gamma} \dar{\psi_2} & L_{0} \rar \dar{\psi_3} & 0 \\
0 \rar & E \rar{\xi} \dar & F \rar{\zeta} \dar & L \rar \dar & 0 \\
 & 0 & 0 & 0 & 
\end{tikzcd}
\end{equation}
In order to prove its commutativity it is sufficient to show that the following squares are commutative 
\begin{equation}
\begin{tikzcd}
E_{0} \rar{\eta} \dar{\phi_1 \circ \psi_1} & F_{0} \dar{\phi_2 \circ \psi_2} \\
E_{0} \rar{\eta} & F_{0} 
\end{tikzcd}
\ \ \ \ \ \ \ \ \ \ \ \ \ \ \
\begin{tikzcd}
F_{0} \rar{\gamma} \dar{\phi_2 \circ \psi_2} & L_{0} \dar{\phi_3 \circ \psi_3} \\
F_{0} \rar{\gamma} & L_{0} 
\end{tikzcd}
\end{equation}
The commutativity of these squares can be derived from the commutative diagram (\ref{diagram_over_dual}) as follows
\begin{equation}
\eta \circ \phi_1 \circ \psi_1 = \phi_2 \circ \xi \circ \psi_1 = \phi_2 \circ \psi_2 \circ \eta,
\end{equation}
\begin{equation}
\gamma \circ \phi_2 \circ \psi_2 = \phi_3 \circ \zeta \circ \psi_1 = \phi_3 \circ \psi_3 \circ \gamma.
\end{equation}
\hfill$\Box$

\noindent Finally, using the previous lemma we are able to prove the main result of the paper.

\begin{Theorem}\label{Non-reduced}
For any sheaf $[E] \in \mathcal{C}$ the following Yoneda product
$$
\cup: \emph{Ext}^{1}(E, E) \times \emph{Ext}^{1}(E, E) \longrightarrow \emph{Ext}^{2}(E, E)
$$
is non-trivial. Therefore, the first obstruction map $\mathcal{Y}_{E, 2}$ is non-zero.
\end{Theorem}
\noindent\textit{Proof:}
By the definition of the subset $\mathcal{C}$, a sheaf $[E] \in \mathcal{C}$ fits into the exact sequence
$$
0 \longrightarrow E \longrightarrow 2 \mathcal{O}_{\mathbb{P}^{3}} \overset{\phi}{\longrightarrow} L \longrightarrow 0,
$$
where $[L] \in \mathcal{P}$ is a line bundle over a smooth curve $C$ from $\mathcal{H}$. Using Lemma \ref{Mumford_obstruction} we can fix the deformation $\textbf{L} \rightarrow L$ of the sheaf $L$ over dual numbers such that the Yoneda square of the corresponding element of the group $\text{Ext}^{1}(L, L)$ is non-zero. Consider also the trivial deformation $2\mathcal{O}_{\mathbb{P}^{3}} \boxtimes \mathcal{O}_{D} \rightarrow 2\mathcal{O}_{\mathbb{P}^{3}}$ and arbitrary surjective morphism $\mathfrak{f}: 2\mathcal{O}_{\mathbb{P}^{3}} \boxtimes \mathcal{O}_{D} \rightarrow \textbf{L}$ such that $\mathfrak{f} \otimes k = \phi$. Then the kernel $\textbf{E} := \text{ker}~\mathfrak{f}$ satisfies the following commutative diagram
\begin{equation}\label{main_triple_over_dual}
\begin{tikzcd}[row sep=scriptsize, column sep=scriptsize]
0 \rar & \textbf{E} \arrow[r] \arrow[d] & 2 \mathcal{O}_{\mathbb{P}^{3}} \boxtimes \mathcal{O}_{D} \arrow[r, "\mathfrak{f}"] \arrow[d] & \textbf{L}  \arrow[r] \arrow[d] & 0 \\
0 \rar & E \arrow[r] & 2 \mathcal{O}_{\mathbb{P}^{3}} \arrow[r, "\phi"] & L \arrow[r] & 0
\end{tikzcd}
\end{equation}
Since $2\mathcal{O}_{\mathbb{P}^{3}} \boxtimes \mathcal{O}_{D}$ and $\textbf{L}$ are flat over $D$, then $\textbf{E}$ is also flat over $D$. Moreover, the equality $\mathfrak{f} \otimes k = \phi$ implies that the morphism $\textbf{E} \rightarrow E$ induces an isomorphism $\textbf{E} \otimes k \simeq E$. Therefore, $\textbf{E} \rightarrow E$ is the deformation of the sheaf $E$ and it defines an element of the group $\text{Ext}^{1}(E, E)$.

Similar to (\ref{extensions}) consider the corresponding exact triples for the sheaves $\text{pr}_{*}\textbf{E}$, $\text{pr}_{*}(2 \mathcal{O}_{\mathbb{P}^{3}} \boxtimes \mathcal{O}_{D})$, $\text{pr}_{*}\textbf{L}$ as 3-term complexes and denote them by $\mathcal{E}, \mathcal{F}, \mathcal{L}$, respectively:
$$
\mathcal{E} \ : \ 0 \longrightarrow E \longrightarrow \text{pr}_{*}\textbf{E} \longrightarrow E \longrightarrow 0,
$$
$$
\mathcal{F} \ : \ 0 \longrightarrow 2 \mathcal{O}_{\mathbb{P}^{3}} \longrightarrow \text{pr}_{*} \big( 2\mathcal{O}_{\mathbb{P}^{3}} \boxtimes \mathcal{O}_{D} \big) \longrightarrow 2 \mathcal{O}_{\mathbb{P}^{3}} \longrightarrow 0,
$$
$$
\mathcal{L} \ : \ 0 \longrightarrow L \longrightarrow \text{pr}_{*}\textbf{L} \longrightarrow L \longrightarrow 0.
$$
Note that due to our assumption we have that
\begin{equation}\label{assumption}
[\mathcal{L} \cup \mathcal{L}] \neq 0 \in \text{Ext}^{2}(L, L).
\end{equation}
Taking pushforward of the upper exact triple of (\ref{main_triple_over_dual}) we obtain the following commutative diagram
\begin{equation}
\begin{tikzcd}[row sep=scriptsize, column sep=scriptsize]
0 \rar & E \arrow[r] \arrow[d] & 2 \mathcal{O}_{\mathbb{P}^{3}} \arrow[r, "\phi"] \arrow[d] & L  \arrow[r] \arrow[d] & 0 \\
0 \rar & \text{pr}_{*}\textbf{E} \arrow[r] \arrow[d] & \text{pr}_{*} \big( 2\mathcal{O}_{\mathbb{P}^{3}} \boxtimes \mathcal{O}_{D} \big) \arrow[r, "\mathfrak{f}"] \arrow[d] & \text{pr}_{*}\textbf{L}  \arrow[r] \arrow[d] & 0 \\
0 \rar & E \arrow[r] & 2 \mathcal{O}_{\mathbb{P}^{3}} \arrow[r, "\phi"] & L \arrow[r] & 0
\end{tikzcd}
\end{equation}
which can be read as the exact triple of 3-complexes
$$
0 \longrightarrow \mathcal{E} \longrightarrow \mathcal{F} \longrightarrow \mathcal{L} \longrightarrow 0.
$$
Next, from Lemma \ref{square_sequence} it follows that the Yoneda squares also fit into the exact triple (\ref{Yoneda_square_formula}), so we have 
$$
0 \longrightarrow \mathcal{E} \cup \mathcal{E} \longrightarrow \mathcal{F} \cup \mathcal{F} \longrightarrow \mathcal{L} \cup \mathcal{L} \longrightarrow 0.
$$

Now suppose that $[\mathcal{E} \cup \mathcal{E}] \in \text{Ext}^{2}(E, E)$ is trivial. By definition it means that there is a commutative diagram (roof) of the following form
\begin{equation}\label{trivialization_E}
\begin{tikzcd}
0 \rar & E \arrow[r, equal] & E \arrow[r, "0"] & E \arrow[r, equal] & E \rar & 0 \\
0 \rar & E \rar \arrow[u, equal] \arrow[d, equal] & E_{0} \rar \arrow[u] \dar & E_{1} \rar \arrow[u] \dar & E \rar \arrow[u, equal] \arrow[d, equal] & 0 \\
0 \rar & E \rar & \text{pr}_{*}\textbf{E} \rar & \text{pr}_{*}\textbf{E} \rar & E \rar & 0 
\end{tikzcd}
\end{equation}
Denote the middle cochain complex of this diagram by $\widetilde{\mathcal{E}}$. Taking the composition of morphisms $\mathcal{E} \rightarrow \mathcal{F}$ and $\widetilde{\mathcal{E}} \rightarrow \mathcal{E}$ we obtain the morphism $\widetilde{\mathcal{E}} \rightarrow \mathcal{F}$. This morphism and its cokernel complex $\widetilde{\mathcal{L}} := \text{coker}(\widetilde{\mathcal{E}} \rightarrow \mathcal{F})$ can be included using the Snake lemma into the following commutative diagram
\begin{equation}\label{aux}
\begin{tikzcd}
& \widetilde{\mathcal{E}} \rar \dar & \mathcal{F} \cup \mathcal{F} \rar \dar[equal] & \widetilde{\mathcal{L}} \rar \dar & 0 \\
0 \rar & \mathcal{E} \cup \mathcal{E} \rar & \mathcal{F} \cup \mathcal{F} \rar & \mathcal{L} \cup \mathcal{L} \rar & 0
\end{tikzcd}
\end{equation}
Assume that the 4-complex $\widetilde{\mathcal{L}}$ has the following form
\begin{equation}\label{L_complex}
\widetilde{\mathcal{L}} \ \ : \ \ 0 \longrightarrow L \longrightarrow L_{0} \longrightarrow L_{1} \overset{\eta}{\longrightarrow} L \longrightarrow 0.
\end{equation}
Note that from the commutative square of the commutative diagram (\ref{trivialization_E})
\begin{equation}
\begin{tikzcd}
E \arrow[r, equal] & E  \\
E \rar \arrow[u, equal] & E_{0} \arrow[u, "\phi"]
\end{tikzcd}
\end{equation}
follows existence of a splitting morphism $\phi : E_{0} \twoheadrightarrow E$, so $E_{0} \simeq E \oplus \widetilde{E}$ where $\widetilde{E}$ is some sheaf. Consider the upper exact triple of the diagram (\ref{aux}) as the following diagram of sheaves
\begin{equation}
\begin{tikzcd}
0 \rar & E \rar \dar & 2 \mathcal{O}_{\mathbb{P}^{3}} \rar \dar & L \rar \dar & 0 \\
& E_{0} \rar \arrow[u, dotted, bend right=50, "\phi"] \dar & 2 \mathcal{O}_{\mathbb{P}^{3}} \oplus 2 \mathcal{O}_{\mathbb{P}^{3}} \rar \arrow[u, dotted, bend right=50, "\psi"] \dar & L_{0} \dar \rar & 0 \\
& E_{1} \rar \dar & 2 \mathcal{O}_{\mathbb{P}^{3}} \oplus 2 \mathcal{O}_{\mathbb{P}^{3}} \rar \dar & L_{1} \rar \dar & 0 \\
0 \rar & E \rar & 2 \mathcal{O}_{\mathbb{P}^{3}}  \rar & L \rar & 0
\end{tikzcd}
\end{equation}
Again using Snake lemma the splitting morphisms $\phi$ and $\psi$ imply the existence of a splitting morphism $\zeta : L_{0} \twoheadrightarrow L$, so we have an isomorphism $L_{0} \simeq L \oplus \widetilde{L}$ for some sheaf $\widetilde{L}$. The morphism $\zeta$ can be included into the following commutative diagram
\begin{equation}
\begin{tikzcd}
0 \rar & L \rar \arrow[d, equal] & L_{0} \rar \arrow[d, "\zeta"] \dar & L_{1} \rar["\eta"] \arrow[d, "\eta"] & L \rar \arrow[d, equal] & 0 \\
0 \rar & L \arrow[r, equal] & L \arrow[r, "0"] & L \arrow[r, equal] & L \rar & 0 
\end{tikzcd}
\end{equation}
where the morphism $\eta$ comes from the definition (\ref{L_complex}) of the complex $\widetilde{\mathcal{L}}$. This means that the equivalence class $[\widetilde{\mathcal{L}}]$ in $\text{Ext}^{2}(L, L)$ is zero. Since the extensions $\mathcal{L}$ and $\widetilde{\mathcal{L}}$ are equivalent, the class $[\mathcal{L} \cup \mathcal{L}] \in \text{Ext}^{2}(L, L)$ is also zero. However, it contradicts to the assumption (\ref{assumption}), so the class $[\mathcal{E} \cup \mathcal{E}] \in \text{Ext}^{2}(E, E)$ cannot be trivial.
\hfill$\Box$

\begin{Corollary}\label{Non-reduced}
Let $\mathfrak{V} \subset \mathcal{M}(14)$ be a maximal closed subscheme of $\mathcal{M}(14)$ containing the closed subset $\overline{\mathcal{C}} \subset \mathcal{M}(14)$. Then the underlying topological subspace of $\mathfrak{V}$ coincides with $\overline{\mathcal{C}}$. Moreover, $\mathfrak{V}$ is generically non-reduced.
\end{Corollary}
\noindent\textit{Proof:} Consider the obstruction map $\mathcal{Y} = \sum\limits_{i=1}^{\infty} \mathcal{Y}_{i}$ defined in Lemma \ref{obstruction_map}. Theorem \ref{Non-reduced} states that $\mathcal{Y}_{2} \neq 0$, so $\mathcal{Y} \neq 0$ and $\text{dim}~\mathcal{Y}^{-1}(0) < 133$. On the other hand, from Lemma \ref{obstruction_map} it follows that the germ of scheme $(\mathcal{Y}^{-1}(0), 0)$ can be immersed in the germ $(\mathcal{M}(14), [E])$. However, the germ $(\mathcal{C}, [E])$ is subgerm of $(\mathcal{M}(14), [E])$, so we have the following inequality 
$$
132 = \text{dim}~\mathcal{C} \leq \text{dim}~\mathcal{Y}^{-1}(0) < 133.
$$
Therefore,we obtain that $\text{dim}~\mathcal{Y}^{-1}(0) = 132$. This means that the subset $\mathcal{C}$ is not contained in some component of the scheme $\mathcal{M}(14)$ whose dimension is more than 132. So the underlying topological space of $\mathfrak{V}$ coincides with $\overline{\mathcal{C}}$. Also as it was shown in Theorem \ref{Main result}, we have $\text{dim}~T_{[E]}\mathcal{M}(14) = 133$ for any sheaf $[E] \in \mathcal{C}$. So the subscheme $\mathfrak{V} \subset \mathcal{M}(14)$ is an irreducible component which is generically non-reduced.
\hfill$\Box$

In conclusion we conjecture that all computations above can be generalized such that the constructed component will be included into a series of generically non-reduced moduli components. More precisely, consider some infinite series of components $\{ \mathcal{R}_{i} \}$ of moduli scheme of reflexive sheaves and series (finite or infinite) of generically non-reduced components $\{ \mathcal{H}_{j}\}$ of the Hilbert scheme which parameterize smooth curves in $\mathbb{P}^{3}$. Now we can construct a family of stable sheaves by taking elementary transforms of the form 
\begin{equation}
0 \longrightarrow E \longrightarrow F \longrightarrow L \oplus \mathcal{O}_{W} \longrightarrow 0,
\end{equation}
where $[F] \in \mathcal{R}_{i}$, $L$ is a line bundle over curve $C \in \mathcal{H}_{j}$ and $W$ is a 0-dimensional subscheme of $\mathbb{P}^{3}$. If we now provide the second condition from (\ref{main_curve_eqs}) and some other mild conditions (see \cite[Lemma 4]{I}) then the isomorphism (\ref{main_isom}) will hold. This means that the difference between dimension of the tangent space $T_{[E]} \mathcal{M}$ and dimension of the family itself is equal to $\text{dim }T_{C}\mathcal{H}_{j} - \text{dim }(\mathcal{H}_{j})_{red}$. So there should exist obstructed deformations of the stable sheaf $E$. Therefore, the corresponding component of the Gieseker-Maruyama moduli scheme is also generically non-reduced.

\end{document}